\documentclass[12pt]{amsart}

\usepackage{amsmath}
\usepackage{amsxtra}
\usepackage{amscd}
\usepackage{amsthm}
\usepackage{amsfonts}
\usepackage{amssymb}
\usepackage{eucal}
\usepackage[matrix,arrow,curve]{xy}
\usepackage[dvips]{graphicx}
\usepackage[dvips]{graphics}
\usepackage[T2A]{fontenc}
\usepackage[cp866]{inputenc}

\newcommand{\tr}{\mathop{\sf tr}\nolimits}
\newcommand{\str}{\mathop{\sf str}\nolimits}
\newcommand{\End}{\mathop{\sf End}\nolimits}

\newcommand{\Hom}{\mathop{\sf Hom}\nolimits}
\newcommand{\Res}{\mathop{\sf Res}\nolimits}

\theoremstyle{plain}
\newtheorem{Thm}[subsection]{Theorem}
\newtheorem{Cor}[subsection]{Corollary}
\newtheorem{Lem}[subsection]{Lemma}
\newtheorem{Prop}[subsection]{Proposition}
\newtheorem{Conj}[subsection]{Conjecture}
\newtheorem{Ex}[subsection]{Example}

\theoremstyle{definition}
\newtheorem{Def}[subsection]{Definition}

\theoremstyle{remark}

\newtheorem{Rem}[subsection]{Remark}

\newtheorem{conj}{Conjecture}

\numberwithin{equation}{section}

\newif\ifShowLabels
\ShowLabelstrue
\newdimen\theight
\def\TeXref#1{%
    \leavevmode\vadjust{\setbox0=\hbox{{\tt
        \quad\quad  {\small \rm #1}}}%
    \theight=\ht0
    \advance\theight by \lineskip
    \kern -\theight \vbox to
    \theight{\rightline{\rlap{\box0}}%
    \vss}%
    }}%

\ShowLabelsfalse

\renewcommand{\sec}[2]{\section{#2}\label{S:#1}%
    \ifShowLabels \TeXref{{S:#1}} \fi}
\newcommand{\ssec}[2]{\subsection{#2}\label{SS:#1}%
    \ifShowLabels \TeXref{{SS:#1}} \fi}

\newcommand{\refs}[1]{Section ~\ref{S:#1}}

\newcommand{\reft}[1]{Theorem ~\ref{T:#1}}
\newcommand{\refl}[1]{Lemma ~\ref{L:#1}}
\newcommand{\refp}[1]{Proposition ~\ref{P:#1}}

\newcommand{\refd}[1]{Definition ~\ref{D:#1}}

\newcommand{\refe}[1]{\eqref{E:#1}}

\newcommand{\refex}[1]{Example ~\ref{Exx:#1}}

\newenvironment{thm}[1]%
    { \begin{Thm} \label{T:#1}  \ifShowLabels \TeXref{T:#1} \fi }%
    { \end{Thm} }

\renewcommand{\th}[1]{\begin{thm}{#1} \sl }
\renewcommand{\eth}{\end{thm} }

\newenvironment{lemma}[1]%
    { \begin{Lem} \label{L:#1}  \ifShowLabels \TeXref{L:#1} \fi }%
    { \end{Lem} }
\newcommand{\lem}[1]{\begin{lemma}{#1} \sl}
\newcommand{\elem}{\end{lemma}}

\newenvironment{propos}[1]%
    { \begin{Prop} \label{P:#1}  \ifShowLabels \TeXref{P:#1} \fi }%
    { \end{Prop} }
\newcommand{\prop}[1]{\begin{propos}{#1}\sl }
\newcommand{\eprop}{\end{propos}}

\newenvironment{corol}[1]%
    { \begin{Cor} \label{C:#1}  \ifShowLabels \TeXref{C:#1} \fi }%
    { \end{Cor} }
\newcommand{\cor}[1]{\begin{corol}{#1} \sl }
\newcommand{\ecor}{\end{corol}}

\newenvironment{defeni}[1]%
    { \begin{Def} \label{D:#1}  \ifShowLabels \TeXref{D:#1} \fi }%
    { \end{Def} }
\newcommand{\defe}[1]{\begin{defeni}{#1} \sl }
\newcommand{\edefe}{\end{defeni}}

\newenvironment{remark}[1]%
    { \begin{Rem} \label{R:#1}  \ifShowLabels \TeXref{R:#1} \fi }%
    { \end{Rem} }
\newcommand{\rem}[1]{\begin{remark}{#1}}
\newcommand{\erem}{\end{remark}}

\newenvironment{conjec}[1]%
    { \begin{Conj} \label{Co:#1}  \ifShowLabels \TeXref{Co:#1} \fi }%
    { \end{Conj} }
\renewcommand{\conj}[1]{\begin{conjec}{#1} \sl }
\newcommand{\econj}{\end{conjec}}

\newenvironment{example}[1]%
    { \begin{Ex} \label{Exx:#1}  \ifShowLabels \TeXref{Exx:#1} \fi }%
    { \end{Ex} }
\newcommand{\ex}[1]{\begin{example}{#1} \sl }
\newcommand{\eex}{\end{example}}

\newcommand{\eq}[1]%
    { \ifShowLabels \TeXref{E:#1} \fi
       \begin{equation} \label{E:#1} }
\newcommand{\eeq}{ \end{equation} }

\newcommand{\prf}{ \begin{proof} }
\newcommand{\epr}{ \end{proof} }

\newcommand\del{\delta}

\newcommand\tet{\theta}

\newcommand\lam{\lambda}        \newcommand\Lam{\Lambda}

\newcommand\calW{{\mathcal{W}}}

 \newcommand\grg{{\mathfrak{g}}}

\newcommand\sdp{\times \hskip -0.3em {\raise 0.3ex
\hbox{$\scriptscriptstyle |$}}} 

\newcommand\Gr{\operatorname{Gr}}

\newcommand\Int{\operatorname{Int}}
\newcommand\Mat{\operatorname{Mat}}

\newcommand\x{\times}
\newcommand\ten{\otimes}

\newcommand{\ra}{\rangle}
\newcommand{\la}{\langle}

\newcommand\nc{\newcommand}

\newcommand{\iso}{{\stackrel{\sim}{\longrightarrow}}}

\nc\aff{\operatorname{aff}}
\nc\oGr{\overline{\Gr}}
\nc\Bun{\operatorname{Bun}}
\nc\hgrg{\widehat{\grg}}
\renewcommand\Int{\operatorname{Int}}
\nc\bInt{\overline{\Int}}
\nc\hatLam{\widehat{\Lam}}
\nc\bmu{\overline{\mu}}
\nc\bnu{\overline{\nu}}
\nc\blambda{\overline{\lam}}

\nc\ocalW{\overline{\calW}}
\nc\pos{\operatorname{pos}}
\nc\IH{\operatorname{IH}}
\nc\Rep{\operatorname{Rep}}
\nc\Gal{\operatorname{Gal}}
\nc{\tilGr}{\widetilde{\Gr}}

\nc\Pic{\operatorname{Pic}}
\nc\pa{\partial}
\nc\Na{\nabla}

\nc{\HC}{{\mathcal{HC}}}
\nc{\on}{\operatorname}
\nc{\BA}{{\mathbb{A}}}
\nc{\BC}{{\mathbb{C}}}
\nc{\BG}{{\mathbb{G}}}
\nc{\BM}{{\mathbb{M}}}
\nc{\BN}{{\mathbb{N}}}
\nc{\BQ}{{\mathbb{Q}}}
\nc{\BP}{{\mathbb{P}}}
\nc{\BR}{{\mathbb{R}}}
\nc{\BZ}{{\mathbb{Z}}}
\nc{\BS}{{\mathbb{S}}}

\nc{\CA}{{\mathcal{A}}}
\nc{\CB}{{\mathcal{B}}}
\nc{\CalC}{{\mathcal C}}
\nc{\CalD}{{\mathcal D}}
\nc{\CE}{{\mathcal{E}}}
\nc{\CF}{{\mathcal{F}}}
\nc{\CG}{{\mathcal{G}}}
\nc{\CH}{{\mathcal{H}}}
\nc{\CK}{{\mathcal{K}}}
\nc{\CL}{{\mathcal{L}}}
\nc{\CM}{{\mathcal{M}}}
\nc{\CMM}{{\mathcal{M}^{\operatorname{gen}}_\hbar(-\rho)}}
\nc{\CN}{{\mathcal{N}}}
\nc{\CO}{{\mathcal{O}}}
\nc{\CP}{{\mathcal{P}}}
\nc{\CQ}{{\mathcal{Q}}}
\nc{\CR}{{\mathcal{R}}}
\nc{\CS}{{\mathcal{S}}}
\nc{\CT}{{\mathcal{T}}}
\nc{\CU}{{\mathcal{U}}}
\nc{\CV}{{\mathcal{V}}}
\nc{\CW}{{\mathcal{W}}}
\nc{\CX}{{\mathcal{X}}}
\nc{\CY}{{\mathcal{Y}}}
\nc{\CZ}{{\mathcal{Z}}}

\nc{\gen}{{\operatorname{gen}}}
\nc{\cM}{{\check{\mathcal M}}{}}
\nc{\csM}{{\check{\mathcal A}}{}}
\nc{\obM}{{\overset{\circ}{\mathbf M}}{}}
\nc{\oCA}{{\overset{\circ}{\mathcal A}}{}}
\nc{\obA}{{\overset{\circ}{\mathbf A}}{}}
\nc{\ooM}{{\overset{\circ}{M}}{}}
\nc{\osM}{{\overset{\circ}{\mathsf M}}{}}
\nc{\vM}{{\overset{\bullet}{\mathcal M}}{}}
\nc{\nM}{{\underset{\bullet}{\mathcal M}}{}}
\nc{\obD}{{\overset{\circ}{\mathbf D}}{}}
\nc{\cp}{{\overset{\circ}{\mathbf p}}{}}
\nc{\ofZ}{{\overset{\circ}{\mathfrak Z}}{}}

\nc{\fa}{{\mathfrak{a}}}
\nc{\fb}{{\mathfrak{b}}}
\nc{\fg}{{\mathfrak{g}}}
\nc{\fgl}{{\mathfrak{gl}}}
\nc{\fh}{{\mathfrak{h}}}
\nc{\fj}{{\mathfrak{j}}}
\nc{\fm}{{\mathfrak{m}}}
\nc{\fn}{{\mathfrak{n}}}
\nc{\fu}{{\mathfrak{u}}}
\nc{\fp}{{\mathfrak{p}}}
\nc{\frr}{{\mathfrak{r}}}
\nc{\fs}{{\mathfrak{s}}}
\nc{\ft}{{\mathfrak{t}}}
\nc{\fT}{{\mathfrak{T}}}
\nc{\ofT}{{\overline{\mathfrak T}}}
\nc{\ofS}{{\overline{\mathfrak S}}}
\nc{\fsl}{{\mathfrak{sl}}}
\nc{\hsl}{{\widehat{\mathfrak{sl}}}}
\nc{\hgl}{{\widehat{\mathfrak{gl}}}}
\nc{\hg}{{\widehat{\mathfrak{g}}}}
\nc{\chg}{{\widehat{\mathfrak{g}}}{}^\vee}
\nc{\hn}{{\widehat{\mathfrak{n}}}}
\nc{\chn}{{\widehat{\mathfrak{n}}}{}^\vee}

\nc{\fA}{{\mathfrak{A}}}
\nc{\fB}{{\mathfrak{B}}}
\nc{\fD}{{\mathfrak{D}}}
\nc{\fE}{{\mathfrak{E}}}
\nc{\fF}{{\mathfrak{F}}}
\nc{\fG}{{\mathfrak{G}}}
\nc{\fI}{{\mathfrak{I}}}
\nc{\fJ}{{\mathfrak{J}}}
\nc{\fK}{{\mathfrak{K}}}
\nc{\fL}{{\mathfrak{L}}}
\nc{\fM}{{\mathfrak{M}}}
\nc{\fN}{{\mathfrak{N}}}
\nc{\frP}{{\mathfrak{P}}}
\nc{\fS}{{\mathfrak S}}
\nc{\fU}{{\mathfrak{U}}}
\nc{\fZ}{{\mathfrak{Z}}}

\nc{\bb}{{\mathbf{b}}}
\nc{\bc}{{\mathbf{c}}}
\nc{\be}{{\mathbf{e}}}
\nc{\bj}{{\mathbf{j}}}
\nc{\bn}{{\mathbf{n}}}
\nc{\bp}{{\mathbf{p}}}
\nc{\bq}{{\mathbf{q}}}
\nc{\bv}{{\mathbf{v}}}
\nc{\bx}{{\mathbf{x}}}
\nc{\by}{{\mathbf{y}}}
\nc{\bw}{{\mathbf{w}}}
\nc{\bA}{{\mathbf{A}}}
\nc{\bB}{{\mathbf{B}}}
\nc{\bC}{{\mathbf{C}}}
\nc{\bK}{{\mathbf{K}}}
\nc{\bD}{{\mathbf{D}}}
\nc{\bH}{{\mathbf{H}}}
\nc{\bM}{{\mathbf{M}}}
\nc{\bN}{{\mathbf{N}}}
\nc{\bS}{{\mathbf{S}}}
\nc{\bT}{{\mathbf{T}}}
\nc{\bV}{{\mathbf{V}}}
\nc{\bW}{{\mathbf{W}}}
\nc{\bX}{{\mathbf{X}}}
\nc{\bP}{{\mathbf{P}}}
\nc{\bZ}{{\mathbf{Z}}}

\nc{\sA}{{\mathsf{A}}}
\nc{\sB}{{\mathsf{B}}}
\nc{\sC}{{\mathsf{C}}}
\nc{\sD}{{\mathsf{D}}}
\nc{\sF}{{\mathsf{F}}}
\nc{\sK}{{\mathsf{K}}}
\nc{\sM}{{\mathsf{M}}}
\nc{\sO}{{\mathsf{O}}}
\nc{\sQ}{{\mathsf{Q}}}
\nc{\sP}{{\mathsf{P}}}
\nc{\sV}{{\mathsf{V}}}
\nc{\sW}{{\mathsf{W}}}
\nc{\sZ}{{\mathsf{Z}}}
\nc{\sfp}{{\mathsf{p}}}
\nc{\sr}{{\mathsf{r}}}
\nc{\sfb}{{\mathsf{b}}}
\nc{\sfc}{{\mathsf{c}}}
\nc{\sd}{{\mathsf{d}}}
\nc{\sg}{{\mathsf{g}}}
\nc{\sfl}{{\mathsf{l}}}

\nc{\BK}{{\bar{K}}}

\nc{\tA}{{\widetilde{\mathbf{A}}}}
\nc{\tB}{{\widetilde{\mathcal{B}}}}
\nc{\tg}{{\widetilde{\mathfrak{g}}}}
\nc{\tG}{{\widetilde{G}}}
\nc{\TM}{{\widetilde{\mathbb{M}}}{}}
\nc{\tO}{{\widetilde{\mathsf{O}}}{}}
\nc{\tU}{{\widetilde{\mathfrak{U}}}{}}
\nc{\TZ}{{\tilde{Z}}}
\nc{\tZ}{\widetilde{Z}{}}
\nc{\tx}{{\tilde{x}}}
\nc{\tbv}{{\tilde{\bv}}}
\nc{\tfP}{{\widetilde{\mathfrak{P}}}{}}
\nc{\tz}{{\tilde{\zeta}}}
\nc{\tmu}{{\tilde{\mu}}}

\nc{\td}{\ddot{\underline{d}}{}}
\nc{\tzeta}{\widetilde{\zeta}{}}
\nc{\hd}{{\widehat{\underline{d}}}}
\nc{\hG}{{\widehat{G}}}
\nc{\hBP}{\widehat{\mathbb P}{}}
\nc{\hQ}{{\widehat{Q}}}
\nc{\hsM}{\widehat{\mathsf M}{}}
\nc{\hfM}{\widehat{\mathfrak M}{}}
\nc{\hCP}{\widehat{\mathcal P}{}}
\nc{\hCR}{\widehat{\mathcal R}{}}
\nc{\hCS}{{\widehat{\mathcal S}}}
\nc{\hfZ}{\widehat{\mathfrak Z}{}}

\nc{\urho}{\underline{\rho}}
\nc{\uB}{\underline{B}}
\nc{\uC}{{\underline{\mathbb{C}}}}
\nc{\ui}{\underline{i}}
\nc{\ofP}{{\overline{\mathfrak{P}}}}

\nc{\hrho}{{\hat{\rho}}}

\nc{\unl}{\underline}
\nc{\ol}{\overline}
\nc{\one}{{\mathbf{1}}}
\nc{\two}{{\mathbf{t}}}

\nc{\Tot}{{\mathop{\operatorname{\rm Tot}}}}
\nc{\Hilb}{{\mathop{\operatorname{\rm Hilb}}}}
\nc{\CHom}{{\mathop{\operatorname{{\mathcal{H}}\it om}}}}
\nc{\defi}{{\mathop{\operatorname{\rm def}}}}
\nc{\length}{{\mathop{\operatorname{\rm length}}}}

\nc{\Cliff}{{\mathsf{Cliff}}}
\nc{\Fl}{{\mathsf{Fl}}}
\nc{\Fib}{{\mathsf{Fib}}}
\nc{\Coh}{{\mathsf{Coh}}}
\nc{\FCoh}{{\mathsf{FCoh}}}

\nc{\reg}{{\text{\rm reg}}}

\nc{\cplus}{{\mathbf{C}_+}}
\nc{\cminus}{{\mathbf{C}_-}}
\nc{\cthree}{{\mathbf{C}_*}}
\nc{\Qbar}{{\bar{Q}}}

\nc{\bh}{{\bar{h}}}
\nc{\bOmega}{{\overline{\Omega}}}
\nc\tGr{\widetilde{\Gr}}

\nc{\seq}[1]{\stackrel{#1}{\sim}}
\nc\ogu{\overline{G/U}}
\nc\chlam{\check{\lam}}

\nc\St{\operatorname{St}}

\nc\uS{\underline{S}}
\nc\QM{\mathcal{QM}}

\nc\zt{\BZ/2\BZ}

\title{Supersymmetric Semisimple Cardy-Frobenius Algebras}
\author{A. Ionov}
\begin{document}

\begin{abstract}

Cardy-Frobenius algebra is the algebraic structure on the space of states in open-closed topological field theory.
 We prove that every  semisimple  super Cardy-Frobenius algebras is the direct sum of the super Cardy-Frobenius algebras of three simple types. We also apply our results to singularity theory via Landau-Ginzburg models and matrix factorizations.

\end{abstract}

\maketitle




\section{Introduction}

Topological field theories reflect the topological aspects of string theory. It was first introduced by Atiyah (\cite{A}). The corresponding algebraical object for 2D-topologiacl field theories is Frobenius algebras for closed field theories (i.e. those where the worldsheet is the closed Riemann surface) and Cardy-Frobenius algebras for open-closed field theories (i.e. those where the worldsheet is the Riemann surface  with boundary). 
The latter ones were developed and defined in \cite{L} and \cite{M}. The proof of equivalence of the corresponding categories could be found in \cite{AN}. The notion of super Cardy-Frobenius algebra were introduced in \cite{L} in the case of super field theories. 

One of the important examples of the open-closed topological super field theories is provided by the so-called Landau-Ginzburg  models. The mathematical way to formulate this theory is via matrix factorization. This theory is of the great interest  because it is strongly connected with some questions in algebraic  geometry and singularities theory. (See, for example, \cite{O} for these connections.)

In this work, we provide the classification of the semisimple super Cardy-Frobenius algebras. More precisely, the semisimple super Cardy-Frobenius algebra could be decomposed in the sum of the Cardy-Frobenius algebras of three types which we call {\itshape elementary}. This is the super analog of the classification of \cite{AN}.


We also use our results in singularity theory via Landau-Ginzburg models. We are interested in the conjecture of Hailong Dao (\cite{D}, Conjecture 3.15). It was first proved in \cite{MPSW}. Another proof with the use of Landau-Ginzburg models was given in \cite{PV}. Using our results, we provide a simple proof for this conjecture in semisimple case.

In \refs{def} we give the definitions and examples.

In \refs{res} we prove our main classificational result (\reft{cf}). 

In \refs{app} we introduce Landau-Ginzburg theory and apply the results of \refs{res} to singularity theory.

We will use the following conventions and notations:

\ssec{}{Conventions and notations} 

Let $K$ be our basic field. We will assume it to be algebraically closed.

We will denote by $E$ the identity matrix and by $E_{i,j}$ the matrix with the entry $1$ in the cell $(i,j)$ and zeroes in other cells. 

\ssec{}{Acknowledgements} The author would like to thank S.Natanzon for  useful discussions. The author was supported by the grants NSh-5138.2014.1 and RFBR 15-01-09242.

\sec{def}{Super Cardy-Frobenius algebra}

Let us start with the definition:

\defe{CF} (\cite{L})
The {\itshape super Cardy-Frobenius algebra} consists of the following data and conditions:

1) a commutative (in the super sense) finite-dimensional  unital $\zt$-graded algebra $A$ (which we will call a {\itshape bulk algebra}) and a finite-dimensional  unital $\zt$-graded algebra $B$  (which we will call a {\itshape boundary algebra}); we  denote by $|\cdot|\in\{0,1\}$ the degree function for elements of $A$ and $B$ of pure degree;

2) a pair of linear forms on them: $\tet_A\colon A\to K$ and $\tet_B\colon B\to K$ such that billinear forms $\la\cdot,\cdot\ra_A$, and $\la\cdot,\cdot\ra_B$ obtained by relation 
$$\la X,Y\ra\colon=\tet(XY).$$ are nondegenerate and symmetric: \eq{s}\tet(XY)=(-1)^{|X||Y|}\tet(YX),\eeq for elements $X,Y$ of pure degree;

3) a  homomorphism of unital algebras $\tau_*\colon A\to B$, which maps $A$ to the center of $B$ and preserves grading ({\itshape bulk-boundary map});

4) a map of vector spaces $\tau^*\colon B\to A$, conjugated to $\tau_*$ with respect to bilinear forms, which either preserves or reverses grading ({\itshape boundary-bulk map});

5) {\itshape (Cardy condition.)}  $\la\tau^*(X),\tau^*(Y)\ra_A=\str_K(m_{X,Y}),$ where $m_{X,Y}\in\End_K(B)$ is the multiplication map:
\eq{mul}m_{X,Y}:B\to B:f\mapsto(-1)^{|X||Y|+|X||f|}Y\cdot f\cdot X,\eeq for any $X, Y\in B$.
\edefe

\ex{triv}
Algebras $A=K, B=0$ with $\tet_A(1)=\lam\ne0$ form a (super) Cardy-Frobenius algebra
\eex
\ex{Mat}
Let $A=K, B=\Mat(n|m)$, i.e. the algebra of endomorphisms of the graded vector space $K^{n|m}$. Consider the forms and the maps: $$\tet_A(x)=\lam x,\enskip \tet_B(X)=\mu\str(X), \enskip\text{where} \enskip\lam,\mu\in K^{\times};$$ $$\tau_*(x)=xE,\enskip\enskip\tau^*(X)=\frac{\mu}{\lam}\str(X) $$ such that the following relation holds: \eq{c1}\lam=\mu^2.\eeq
Then it is a Cardy-Frobenius algebra.

\eex

\prf
Let us check the Cardy condition. 

We will abstractly identify $\Mat(n|m)$ with $\Mat(n+m)$ as algebra by forgeting grading and consider the elements $E_{i,j}$ as the elements of $\Mat(n|m)$.

By linearity it suffices to check the Cardy condition for elements $E_{i,j}$ and compute the traces in basis  $E_{i,j}$.
Note that: $$E_{i_1,j_1}E_{k,l}E_{i_2,j_2}=\del_{j_1,k}\del_{l,i_2}E_{i_1,j_2}.$$ Thus,
$$\str(m_{E_{i_1,j_1},E_{i_2,j_2}})=\begin{cases}
0, &\text{ if $i_1\ne j_1$ or $i_2\ne j_2$;}\\
1, &\text{if  $i_1=j_1$, $i_2=j_2$ and $|E_{i_1,i_2}|=0$;}\\
-1, &\text{ if $i_1=j_1$, $i_2=j_2$ and $|E_{i_1,i_2}|=1$.}\\

\end{cases}$$

Let us also note that:
$$\tau^*(E_{i,j})=\begin{cases}
\frac{\mu}{\lam}, &\text{if  $i=j\le n$;}\\
-\frac{\mu}{\lam}, &\text{ if $i=j>n$;}\\
0, &\text{ if $i\ne j$.}\end{cases}$$

Then,
$$\la\tau^*(E_{i_1,j_1}),\tau^*(E_{i_2,j_2})\ra=\begin{cases}
0, &\text{ if $i_1\ne j_1$ or $i_2\ne j_2$;}\\
\frac{\mu^2}{\lam}, &\text{if  $i_1=j_1$, $i_2=j_2$ and either $i_1,i_2\le n$, or $i_1,i_2>n$;}\\
-\frac{\mu^2}{\lam}, &\text{otherwise. }\\

\end{cases}$$

Therefore, \refe{c1} is equivalent to the Cardy condition.

\epr
\ex{Q}Let $A=K, B=Q(n):=\Mat(n)[\xi]$, where $|\xi|=1$ and $\xi^2=1$. Consider the forms and the maps:  $$\tet_A(x)=\lam x,\enskip\tet_B(X+Y\xi)=\mu\tr(Y),\enskip\text{where} \enskip\lam,\mu\in K^{\times};$$ $$\tau_*(x)=xE,\enskip\enskip\tau^*(X+Y\xi)=\frac{\mu}{\lam}\tr(Y) $$ such that the following relation holds:\eq{c2}\lam=\frac{1}{2}\mu^2.\eeq
Then it is a Cardy-Frobenius algebra.
\eex
\prf
Let us check the Cardy condition. By linearity it suffices to consider only elements $E_{i,j}$ and $E_{i,j}\xi$ and compute traces in the same basis.

In this basis operator $m_{E_{i_1,j_1},E_{i_2,j_2}}$ has nonzero diagonal entries if and only if $i_1=j_1$ and $i_2=j_2$. In this case, it is $1$ for $E_{i_1, i_2}$ and $1$ for $E_{i_1, i_2}\xi$. Since one of them is even and the other is odd, they are taken with different signs and the right hand side of the Cardy condition vanishes. 

If $X$ and $Y$ have the different parity, then in the chosen basis, operator $m_{X,Y}$ has no nonzero entries on the diagonal.

Note that $\tau^*$ vanishes on even elements. Thus, the Cardy condition is valid if at least one of the elements is even. Therefore, it remains to check it for $X=E_{i_1,j_1}\xi, Y=E_{i_2,j_2}\xi$.

Operator $m_{E_{i_1,j_1}\xi,E_{i_2,j_2}\xi}$ has non zero diagonal entries in the chosen basis only if and only if $i_1=j_1$ and $i_2=j_2$. In this case, it is $1$ for $E_{i_1, i_2}$ and  $1$ for $E_{i_1, i_2}\xi$. Due to the choice of signs in \refe{mul}, both of them are taken with the sign $+$ and the right-hand side of the Cardy condition is equal to $2$.

Let us note that $\tau^*(E_{i_1,i_1}\xi)=\tau^*(E_{i_2,i_2}\xi)=\frac{\mu}{\lam}$ and the left-hand side of the Cardy condition equals $\frac{\mu^2}{\lam}$. Therefore, \refe{c2} is equivalent to the Cardy condition.

\epr


 

\defe{} We will call the Cardy-Frobenius algebras of \refex{triv}, \refex{Mat}, \refex{Q} elementary Cardy-Frobenius algebras.
\edefe

Further, we will work with {\itshape semisimple super Cardy-Frobenius algebras}:
\defe{}
The Cardy-Frobenius algebra is called  {\itshape semisimple} if $A$ and $B$ are semisimple in the category of $\zt$-graded algebras.
\edefe
\ex{} The elementary Cardy-Frobenius algebras are semisimple.
\eex
\defe{}
The {\itshape direct sum} of the Cardy-Frobenius algebras $(A_1,B_1,\tet_{A_1},\tet_{B_1},\tau^*_1,\tau_{1*})$ and $(A_2,B_2,\tet_{A_2},\tet_{B_2},\tau^*_2,\tau_{2*})$ is defined by taking $A_1\oplus A_2$ as a bulk algebra, $B_1\oplus B_2$ as a boundary algebra and $$\tet_{A_1\oplus A_2}(a_1\oplus a_2)=\tet_{A_1}(a_1)+\tet_{A_2}(a_2),\tet_{B_1\oplus B_2}(b_1\oplus b_2)=\tet_{B_1}(b_1)+\tet_{B_2}(b_2),$$ $$\tau^*(b_1\oplus b_2)=\tau^*_1(b_1)+\tau^*_2(b_2), \tau_*(a_1\oplus a_2)=\tau_{1*}(a_1)+\tau_{2*}(a_2).$$
\edefe


\sec{res}{Classification of semisimple Cardy-Frobenius algebras}

In this section we provide a full classification of semisimple Cardy-Frobenius algebras. It is the super analog of the classification of \cite{AN}.

\th{cf}
Semisimple Cardy-Frobenius algebra is the direct sum of the elementary Cardy-Frobenius algebras.
\eth

Let $(A,B,\tet_A,\tet_B,\tau^*,\tau_*)$ be a semisimple  Cardy-Frobenius algebra.

\ssec{alg}{Algebras} First of all we need the classification of semisimple $\zt$-graded algebras. We will use the following analog of the Artin-Wedderburn theorem:

\th{AW}(\cite{W}) 

1) Semisimple  $\zt$-graded algebras are the direct sums of simple.

2) Simple  $\zt$-graded algebras are either $\Mat(n|m)$ or $Q(n)$.


\eth

Note that $K=\Mat(1|0)$ is the only commutative algebra among the list of simple algebras. Thus we obtain the following:
\cor{}
We have the isomorphism of the algebras  $A=K\oplus\ldots\oplus K$. In particular, $A$ is purely even.
\ecor

We will denote by $e_1,\ldots,e_k\in A$ the idempotents corresponding to the decomposition $A=K\oplus\ldots\oplus K$.



\ssec{form}{Linear forms} 

Secondly, we need to find linear forms on semisimple algebras satisfying \refe{s} and providing nondegenerate bilinear forms.



According to \reft{AW},  it is sufficient to study linear forms satisfying \refe{s} for $\Mat(n|m)$ and $Q(n)$ as we can restrict $\tet$  on direct summands. The answer in this cases is given by \refl{f1} and \refl{f2}. 

\lem{f1} 
Linear form on $\Mat(n|m)$ satisfying \refe{s} is proportional to $\str(\cdot)$.
\elem
\prf
Consider the even part of the algebra. It is isomorphic  to $\Mat(n)\oplus\Mat(m)$. It is well-known that  the form on matrix algebra satisfying $\tet(XY)=\tet(YX)$ is proportional to trace. Then the restriction of $\tet$ on each of summands is proportional to trace. We will asume that $n,m>0$.

Let us prove that these coefficients of proportionality are opposite.


Note that $$E_{1,1}=E_{1,n+1}E_{n+1,1} \enskip \text{and}\enskip E_{n+1,n+1}=E_{n+1,1}E_{1,n+1}.$$ Then the condition \refe{s} implies

$$\tet(E_{1,1})=-\tet(E_{n+1,n+1}).$$ 

Therefore, the restriction of $\tet(\cdot)$ on the even part is proportional to $\str(\cdot)$. 

Thus, it remains to check that the restriction to the odd part vanishes. Let $X\in\Mat(n), Y\in\Mat(m)$ be invertible elements, then
$$\begin{pmatrix}X & \ 0  \\ 0 & \ Y^{-1}  \end{pmatrix}\begin{pmatrix}0 \ Z\\ 0 \ 0 \end{pmatrix}\begin{pmatrix}X^{-1} &\ 0 \\ 0 &\ Y \end{pmatrix}=\begin{pmatrix} 0 &\  XZY  \\ 0 & \ 0 \end{pmatrix}.$$ Then it follows that $$\tet(\begin{pmatrix}0 \ Z\\ 0 \ 0 \end{pmatrix})=\tet(\begin{pmatrix} 0 &\  XZY  \\ 0 & \ 0 \end{pmatrix}).$$ But there are no nonzero forms on $\Mat_{n\x m}$ satisfying this property for every invertible 
$X, Y$.
\epr

\lem{f2}
Linear form on $Q(n)$ satisfying \refe{s} is proportional to the form of \refex{Q}.
\elem
\prf
Note that for $X\in\Mat(n)$: $$\tet(X)=\tet(E\xi\cdot X\xi)=\tet(X\xi\cdot E\xi).$$
The last two expressions are opposite by \refe{s} and, thus, it is zero. Therefore, the restriction of $\tet$ to the even part is identically zero. 

The odd part of $Q(n)$ consists of the elements $X\xi$, where $X\in\Mat(n)$. Thus, we need the linear map $\tet$ from elements of the form $X\xi$ such that $$\tet(Y\cdot Z\xi)=\tet(Z\xi\cdot Y),$$ for $Y,Z\in\Mat(n)$. But this is exactly the same as the linear form on $\Mat(n)$ such that $\tet(YZ)=\tet(ZY)$. All such linear forms are proportional to $\tr(\cdot)$, as was discussed earlier.


\epr

\ssec{map}{The maps $\tau_*$ and  $\tau^*$ and the Cardy condition}

Let $B=B_1\oplus\ldots\oplus B_m$ be the decomposition of $B$ in the sum of simple $\zt$-graded algebras. Denote by $1_i$ the unit of $B_i$ as an element of $B$.
\prop{dec}
There is the reordering of $e_i$ and $B_j$ such that 

1) $\tau_*(e_i)=1_i$ for $i=1,\ldots,m$;

2) $\tau_*(e_i)=0$ for $i=m+1,\ldots,k$.

In particular, $k\ge m$.
\eprop
\prf
By \refd{CF} 4), the images of the elements $e_i$ under $\tau_*$ lie in the center of $B$. The centers of simple algebras $\Mat(n|m)$ and $Q(n)$ are generated by the unit element. Thus, the center of $B$ is $K1_1\oplus\ldots\oplus K1_m$.

Moreover, by \refd{CF} 4), $\tau_*$ is the homomorphism of unital algebras. It follows that $\tau_*(e_i)=1_{i_1}+\ldots+1_{i_l}$ for some $\{i_1,\ldots,i_l\}\subset\{1,\ldots,m\}$. These subsets do not intersect for different $i$, since $e_ie_j=0$ for $i\ne j$.

It remains to check that these subsets for each $i$ have no more than one element. Assume that $\tau_*(e_i)=1_{i_1}+1_{i_2}+\ldots$. Then, for any pair of elements $X\in B_{i_1}, Y\in B_{i_2}$ the right-hand side of the Cardy condition (\refd{CF} 5) ) vanishes, since $YBX\subset B_{i_1}\cap B_{i_2}=\{0\}$. 

On the other hand, by conjugacy $\tau^*(X)=\frac{\la 1_{i_1},X\ra_B}{\la e_i, e_i\ra_A}e_i$ and $\tau^*(Y)=\frac{\la 1_{i_2},Y\ra_B}{\la e_i, e_i\ra_A}e_i$. Thus, the left hand side of the Cardy condition equals $$\frac{\la 1_{i_1},X\ra_B\la 1_{i_2},Y\ra_B}{\la e_i, e_i\ra_A}.$$ But this is nonzero for some $X,Y$, because the form $\la\cdot,\cdot\ra_B$ is nondegenerate.
\epr
\prf[Proof of \reft{cf}]
By \refp{dec} we have a decomposition of the semisimple Cardy-Frobenius algebra into the direct sum of Cardy-Frobenius algebras with $A=K$ and $B$ equal to $0,\Mat(n|m)$ or $Q(n)$  and with $\tau_*$ mapping the unit of $A$ to the unit of $B$.

\refl{f1} and \refl{f2} imply that the linear forms coincide with the forms of \refex{Mat}, \refex{Q}. 

Also note that $\tau^*$ is uniquely defined from $\tau_*$ by conjugacy.

The checking of the Cardy condition for \refex{Mat} and \refex{Q} implies that it is equivalent to \refe{c1} and \refe{c2}.
\epr

\sec{app}{Matrix factorizations}

We are substantially interested in the open-closed topological field theories coming from matrix factorization. We refer to \cite{CM},\cite{O},\cite{PV} for defintions of the category of matrix factorizations and its basic properties.

Let  $R$ be a matrix factorization of the polynomial $W\in K[x_1,\ldots,x_n]$. Then put $A=K[x_1,\ldots,x_n]/(\partial_{x_1}W,\ldots,\partial_{x_n}W)$, i.e. the Milnor ring of $W$ and $B=\End(R)$, i.e. the cohomologies of endomorphisms dg-ring of $R$ in the dg-category of matrix factorizations of $W$. 
Let us also define the maps: $$\tau_*\colon X\mapsto X\cdot1_R ;\enskip \tau^*\colon Y\mapsto(-1)^{\binom{n+1}{2}}\mathrm{str}(Y\partial_{x_1}d_R\ldots\partial_{x_n}d_R)
$$
and linear forms:
$$
\tet_A(X)=\Res_{K[x]/K}\Bigl[\frac{X{\mathrm d}x}{\partial_{x_1}W,\ldots,\partial_{x_n}W}\Bigr] .
$$
and
$$
\tet_B(Y)=\Res_{K[x]/K}\Bigl[\frac{\mathrm{str}(Y\partial_{x_1}d_R\cdots\partial_{x_n}d_R){\mathrm d}x}{\partial_{x_1}W,\ldots,\partial_{x_n}W}\Bigr].
$$
This forms naturally come from the structure of Calabi-Yau dg-category of dimension $n$ on the category of matrix factorizations of $W$. 

\prop{mf}
The above $A,B,\tet_A,\tet_B,\tau^*,\tau_*$ are well-defined and form a Cardy-Frobenius algebra.
\eprop
The proofs of \refp{mf} can be found in \cite{CM} and \cite{PV}. Here we give a proof of Cardy in the special but remarkable case.

For $\zt$-graded vector space $V=V_0\oplus V_1$ let us denote by $\chi(V)$ its Euler characteristic:
$$\chi(V)=\dim V_0-\dim V_1.$$

Note that $m_{1,1}$ is the identity operator and so $\str(m_{1,1})$ is equal to the Euler characteristic $\chi(\End(R))$. Let us also note that if $n$ is odd, then  $\tau^*(1)=\mathrm{str}(\partial_{x_1}d_R\ldots\partial_{x_n}d_R)=0$, since this is the product of the odd number of odd elements. Thus, the following result is equivalent to the Cardy condition for $X=Y=1$:

\th{D}
Let $R$ be a matrix factorization of $W\in K[x_1,\ldots,x_n]$ and let $n$ be odd. Then $$\chi(\End(R))=0.$$
\eth

This result is notable because it is equivalent to the conjecture of Hailong Dao (\cite{D}, Conjecture 3.15). Indeed, in \cite{O} the relation between the category of singularities and the category of matrix factorizations is established and the conjecture of \cite{D} could be seen as the vanishing of the Euler characteristics of $\Hom$-complexes in the category of singularities.




We provide a simple proof of this result under an assumption of semisimplicity. By this we means that we will suppose that $B=\End(R)$ is semisimple as the $\zt$-graded algebra.

\prf 
Consider the decomposition of $B$ into the sum of simple algebras: $B=B_1\oplus\ldots\oplus B_n$. Note that, $\tet_B(Y)$ vanishes for every even $Y\in B$ if $n$ is odd. Then, it follows from \refl{f1} and \refl{f2} that $B_i=Q(n_i)$ for all $i$.

Then we have:
$$\chi(B)=\sum\chi(Q(n_i))=0.$$
\epr



\bigskip
\footnotesize{
{\bf A.I.}: National Research University
Higher School of Economics\\
Department of Mathematics, 20 Myasnitskaya st, Moscow 101000 Russia;\\
{\tt 8916456@rambler.ru}}

\end{document}

The problem of this argument is that it is not obvious that the bilinear form $\la\cdot,\cdot\ra_B$ provided by $\tet_B$ is nondegenerate. For this purpose we need to remind another definition of $\tet_B$.

The standard argument (see \cite{PV}) endows the category of matrix factorizations with the structure of Calabi-Yau category of dimension $n$. In particular, it means that there is an isomorphism
$$\End(R)^{\vee}\iso\End(R)[n]$$ or equivalently nondegenerate graded symmetric invariant billinear form 
$$\End(R)\ten\End(R)\to K[n],$$ which gives a linear form $$\tet_B^\prime\colon\End(R)\to K[n].$$

Consider the decomposition of $B$ into the sum of simple algebras: $B=B_1\oplus\ldots\oplus B_n$. Note that, $\tet_B(Y)$ vanishes for every even $Y\in B$ if $n$ is odd. Then, it follows from \refl{f1} and \refl{f2} that $B_i=Q(n_i)$ for all $i$.

Then we have:
$$\chi(B)=\sum\chi(Q(n_i))=0.$$
\epr

\rem{} We introduce $\tet^\prime_B$ here in order to make the proof more self-contained, because, generally, it is not obvious that the form $\tet_B$ defined earlier provides nondegenerate symmetric billinear form, but we do not need to know it for this exact formula.

In fact, $\tet^\prime_B=\tet_B$. The proof of this fact could be found in \cite{PV}.

.............................................................................................................

Again as in non $\zt$-graded situation (see \cite{AN}) it implies that idempotents of bulk algebra maps to not more than one component  of boundary algebrs (as its unit). 

For the $\Mat(n|m)$ case the Cardy condition in terms of metrics has the form \eq{c1}\lam=\mu^2,\eeq where $\lam$ is the length of idempotent and the metric on this component of boundary algebra is $\mu\str(\cdot)$. It follow by consideration of matrices $E_{ij}$. 
Nontrivial identity equivalent to \refe{c1} appears for  $\phi={E_{ii}}, \psi={E_{jj}}$.

For $Q(n)$ let us denote by $\widetilde{E_{ij}}$ the matrix $\begin{pmatrix}0 \ E_{ij}\\ E_{ij} \ 0 \end{pmatrix}$. Let the restriction of metric on odd part of this $Q(n)$ be $\mu\tr(\cdot)$ and let the length of corresponding idempotent be $\lam$. Then the image of $\widetilde{E_{ij}}$ is $\frac{\mu}{\lam}\del_{ij}$ and the Cardy condition for for $\phi=\widetilde{E_{ii}}, \psi=\widetilde{E_{jj}}$ turns to be \eq{c2}\frac{\mu^2}{\lam}=2;\eeq  for  even elements of $Q(n)$ Cardy condition is trivial since they maps to $0$ and the supertrace of multiplication by them vanishes.

.............................................................................................................

\title{Semisimple super Cardy-Frobenius Algebras}
\setcounter{tocdepth}{1}
\begin{document}

\begin{abstract}
Cardy-Frobenius algebra is the algebraic structure on the space of states in open-closed topological field theory.
 We prove that every  semisimple  super Cardy-Frobenius algebras is the direct sum of the Cardy-Frobenius algebras of three simple types. We also apply our results to singularity theory via Landau-Ginzburg models and matrix factorizations.

\end{abstract}

\maketitle






\sec{}{Introduction}

Topological field theories reflect the topological aspects of string theory. It was first introduced by Atiyah (\cite{A}). The corresponding algebraical object for 2D-topologiacl field theories is Frobenius algebras for closed field theories (i.e. those where the worldsheet is the closed Riemann surface) and Cardy-Frobenius algebras for open-closed field theories (i.e. those where the worldsheet is the Riemann surface  with boundary). The latter ones were developed and defined in \cite{L} and \cite{MS}. (See also \cite{C} for more modern, categorical approach.)

One of the important examples of the open-closed topological field theories is provided by the so-called Landau-Ginzburg  models. The mathematical way to formulate this theory is via matrix factorization. This theory is of the great interest  because it is strongly connected with some questions in algebraic  geometry and singularities theory. (See, for example, \cite{O} for these connections.)

In this work, we provide the classification of the semisimple super Cardy-Frobenius algebras. More precisely, the semisimple super Cardy-Frobenius algebra could be decomposed in the sum of the Cardy-Frobenius algebras of three types which we call {\itshape elementary}. This is the super analog of the classification of \cite{AN}.


We also use our results in singularity theory via Landau-Ginzburg models. We are interested in the conjecture of Hailong Dao (\cite{D}, Conjecture 3.15). It was first proved in \cite{MPSW}. Another proof with the use of Landau-Ginzburg models was given in \cite{PV}. Using our results, we provide a simple proof of this conjecture in semisimple case.

In \refs{def} we give the definitions and examples.

In \refs{res} we prove our main classificational result (\reft{cf}). 

In \refs{app} we introduce Landau-Ginzburg theory and apply the results of \refs{res} to singularity theory.

\ssec{}{Conventions and notations} 

Let $K$ be our basic field. We will assume it to be algebraically closed.

We will denote by $E$ the identity matrix and by $E_{i,j}$ the matrix with the entry $1$ in the cell $(i,j)$ and zeroes in other cells.

\sec{def}{Super Cardy-Frobenius algebra}

\defe{CF} (\cite{L})
The {\itshape super Cardy-Frobenius algebra} consists of the following data and conditions:

1) a commutative finite-dimensional  unital $\zt$-graded algebra $A$ (which we will call a {\itshape bulk algebra}) and a finite-dimensional  unital $\zt$-graded algebra $B$  (which we will call a {\itshape boundary algebra}); we  denote by $|\cdot|\in\{0,1\}$ the degree function for elements of $A$ and $B$ of pure degree;

2) a pair of linear forms on them: $\tet_A\colon A\to K$ and $\tet_B\colon B\to K$ such that billinear forms $\la\cdot,\cdot\ra_A$, and $\la\cdot,\cdot\ra_B$ obtained by relation 
$$\la X,Y\ra\colon=\tet(XY).$$ are nondegenerate and symmetric: \eq{s}\tet(XY)=(-1)^{|X||Y|}\tet(YX),\eeq for elements $X,Y$ of pure degree;

3) a  homomorphism of unital algebras $\tau_*\colon A\to B$, which maps $A$ to the center of $B$ and preserves grading ({\itshape bulk-boundary map});

4) a map of vector spaces $\tau^*\colon B\to A$, conjugated to $\tau_*$ with respect to bilinear forms, which either preserves or reverses grading ({\itshape boundary-bulk map});

5) {\itshape (Cardy condition.)}  $\la\tau^*(X),\tau^*(Y)\ra_A=\str_K(m_{X,Y}),$ where $m_{X,Y}\in\End_K(B)$ is the multiplication map:
\eq{mul}m_{X,Y}:B\to B:f\mapsto(-1)^{|X||Y|+|X||f|}Y\cdot f\cdot X,\eeq for any $X, Y\in B$.
\edefe

\ex{triv}
Algebras $A=K, B=0$ with $\tet_A(1)=\lam\ne0$ form a (super) Cardy-Frobenius algebra
\eex
\ex{Mat}
Let $A=K, B=\Mat(n|m)$, i.e. the algebra of endomorphisms of the graded vector space $K^{n|m}$. Consider the forms and the maps: $$\tet_A(x)=\lam x,\enskip \tet_B(X)=\mu\str(X), \enskip\text{where} \enskip\lam,\mu\in K^{\times};$$ $$\tau_*(x)=xE,\enskip\enskip\tau^*(X)=\frac{\mu}{\lam}\str(X) $$ such that the following relation holds: \eq{c1}\lam=\mu^2.\eeq
Then it is a Cardy-Frobenius algebra.

\eex

\prf
Let us check the Cardy condition. 

We will abstractly identify $\Mat(n|m)$ with $\Mat(n+m)$ as algebra by forgeting grading and consider the elements $E_{i,j}$ as the elements of $\Mat(n|m)$.

By linearity it suffices to check the Cardy condition for elements $E_{i,j}$ and compute the traces in basis  $E_{i,j}$.
Note that: $$E_{i_1,j_1}E_{k,l}E_{i_2,j_2}=\del_{j_1,k}\del_{l,i_2}E_{i_1,j_2}.$$ Thus,
$$\str(m_{E_{i_1,j_1},E_{i_2,j_2}})=\begin{cases}
0, &\text{ if $i_1\ne j_1$ or $i_2\ne j_2$;}\\
1, &\text{if  $i_1=j_1$, $i_2=j_2$ and $|E_{i_1,i_2}|=0$;}\\
-1, &\text{ if $i_1=j_1$, $i_2=j_2$ and $|E_{i_1,i_2}|=1$.}\\

\end{cases}$$

Let us also note that:
$$\tau^*(E_{i,j})=\begin{cases}
\frac{\mu}{\lam}, &\text{if  $i=j\le n$;}\\
-\frac{\mu}{\lam}, &\text{ if $i=j>n$;}\\
0, &\text{ if $i\ne j$.}\end{cases}$$

Then,
$$\la\tau^*(E_{i_1,j_1}),\tau^*(E_{i_2,j_2})\ra=\begin{cases}
0, &\text{ if $i_1\ne j_1$ or $i_2\ne j_2$;}\\
\frac{\mu^2}{\lam}, &\text{if  $i_1=j_1$, $i_2=j_2$ and either $i_1,i_2\le n$, or $i_1,i_2>n$;}\\
-\frac{\mu^2}{\lam}, &\text{otherwise. }\\

\end{cases}$$

Therefore, \refe{c1} is equivalent to the Cardy condition.

\epr
\ex{Q}Let $A=K, B=Q(n):=\Mat(n)[\xi]$, where $|\xi|=1$ and $\xi^2=1$. Consider the forms and the maps:  $$\tet_A(x)=\lam x,\enskip\tet_B(X+Y\xi)=\mu\tr(Y),\enskip\text{where} \enskip\lam,\mu\in K^{\times};$$ $$\tau_*(x)=xE,\enskip\enskip\tau^*(X+Y\xi)=\frac{\mu}{\lam}\tr(Y) $$ such that the following relation holds:\eq{c2}\lam=\frac{1}{2}\mu^2.\eeq
Then it is a Cardy-Frobenius algebra.
\eex
\prf
Let us check the Cardy condition. By linearity it suffices to consider only elements $E_{i,j}$ and $E_{i,j}\xi$ and compute traces in the same basis.

In this basis operator $m_{E_{i_1,j_1},E_{i_2,j_2}}$ has nonzero diagonal entries if and only if $i_1=j_1$ and $i_2=j_2$. In this case, it is $1$ for $E_{i_1, i_2}$ and $1$ for $E_{i_1, i_2}\xi$. Since one of them is even and the other is odd, they are taken with different signs and the right hand side of the Cardy condition vanishes. 

If $X$ and $Y$ have the different parity, then in the chosen basis, operator $m_{X,Y}$ has no nonzero entries on the diagonal.

Note that $\tau^*$ vanishes on even elements. Thus, the Cardy condition is valid if at least one of the elements is even. Therefore, it remains to check it for $X=E_{i_1,j_1}\xi, Y=E_{i_2,j_2}\xi$.

Operator $m_{E_{i_1,j_1}\xi,E_{i_2,j_2}\xi}$ has non zero diagonal entries in the chosen basis only if and only if $i_1=j_1$ and $i_2=j_2$. In this case, it is $1$ for $E_{i_1, i_2}$ and  $1$ for $E_{i_1, i_2}\xi$. Due to the choice of signs in \refe{mul}, both of them are taken with the sign $+$ and the right-hand side of the Cardy condition is equal to $2$.

Let us note that $\tau^*(E_{i_1,i_1}\xi)=\tau^*(E_{i_2,i_2}\xi)=\frac{\mu}{\lam}$ and the left-hand side of the Cardy condition equals $\frac{\mu^2}{\lam}$. Therefore, \refe{c2} is equivalent to the Cardy condition.

\epr


 

\defe{} We will call the Cardy-Frobenius algebras of \refex{triv}, \refex{Mat}, \refex{Q} elementary Cardy-Frobenius algebras.
\edefe

Further, we will work with {\itshape semisimple super Cardy-Frobenius algebras}:
\defe{}
The Cardy-Frobenius algebra is called  {\itshape semisimple} if $A$ and $B$ are semisimple in the category of $\zt$-graded algebras.
\edefe
\ex{} The elementary Cardy-Frobenius algebras are semisimple.
\eex
\defe{}
The {\itshape direct sum} of the Cardy-Frobenius algebras $(A_1,B_1,\tet_{A_1},\tet_{B_1},\tau^*_1,\tau_{1*})$ and $(A_2,B_2,\tet_{A_2},\tet_{B_2},\tau^*_2,\tau_{2*})$ is defined by taking $A_1\oplus A_2$ as a bulk algebra, $B_1\oplus B_2$ as a boundary algebra and $$\tet_{A_1\oplus A_2}(a_1\oplus a_2)=\tet_{A_1}(a_1)+\tet_{A_2}(a_2),\tet_{B_1\oplus B_2}(b_1\oplus b_2)=\tet_{B_1}(b_1)+\tet_{B_2}(b_2),$$ $$\tau^*(b_1\oplus b_2)=\tau^*_1(b_1)+\tau^*_2(b_2), \tau_*(a_1\oplus a_2)=\tau_{1*}(a_1)+\tau_{2*}(a_2).$$
\edefe


\sec{res}{Classification of semisimple Cardy-Frobenius algebras}

In this section we provide a full classification of semisimple Cardy-Frobenius algebras. It is the super analog of the classification of \cite{AN}.

\th{cf}
Semisimple Cardy-Frobenius algebra is the direct sum of the elementary Cardy-Frobenius algebras.
\eth

Let $(A,B,\tet_A,\tet_B,\tau^*,\tau_*)$ be a semisimple  Cardy-Frobenius algebra.

\ssec{alg}{Algebras} First of all we need the classification of semisimple $\zt$-graded algebras. We will use the following analog of the Artin-Wedderburn theorem:

\th{AW}(\cite{W}) 

1) Semisimple  $\zt$-graded algebras are the direct sums of simple.

2) Simple  $\zt$-graded algebras are either $\Mat(n|m)$ or $Q(n)$.


\eth

Note that $K=\Mat(1|0)$ is the only commutative algebra among the list of simple algebras. Thus we obtain the following:
\cor{}
We have the isomorphism of the algebras  $A=K\oplus\ldots\oplus K$. In particular, $A$ is purely even.
\ecor

We will denote by $e_1,\ldots,e_k\in A$ the idempotents corresponding to the decomposition $A=K\oplus\ldots\oplus K$.



\ssec{form}{Linear forms} 

Secondly, we need to find linear forms on semisimple algebras satisfying \refe{s} and providing nondegenerate bilinear forms.



According to \reft{AW},  it is sufficient to study linear forms satisfying \refe{s} for $\Mat(n|m)$ and $Q(n)$ as we can restrict $\tet$  on direct summands. The answer in this cases is given by \refl{f1} and \refl{f2}. 

\lem{f1} 
Linear form on $\Mat(n|m)$ satisfying \refe{s} is proportional to $\str(\cdot)$.
\elem
\prf
Consider the even part of the algebra. It is isomorphic  to $\Mat(n)\oplus\Mat(m)$. It is well-known that  the form on matrix algebra satisfying $\tet(XY)=\tet(YX)$ is proportional to trace. Then the restriction of $\tet$ on each of summands is proportional to trace. We will asume that $n,m>0$.

Let us prove that these coefficients of proportionality are opposite.


Note that $$E_{1,1}=E_{1,n+1}E_{n+1,1} \enskip \text{and}\enskip E_{n+1,n+1}=E_{n+1,1}E_{1,n+1}.$$ Then the condition \refe{s} implies

$$\tet(E_{1,1})=-\tet(E_{n+1,n+1}).$$ 

Therefore, the restriction of $\tet(\cdot)$ on the even part is proportional to $\str(\cdot)$. 

Thus, it remains to check that the restriction to the odd part vanishes. Let $X\in\Mat(n), Y\in\Mat(m)$ be invertible elements, then
$$\begin{pmatrix}X & \ 0  \\ 0 & \ Y^{-1}  \end{pmatrix}\begin{pmatrix}0 \ Z\\ 0 \ 0 \end{pmatrix}\begin{pmatrix}X^{-1} &\ 0 \\ 0 &\ Y \end{pmatrix}=\begin{pmatrix} 0 &\  XZY  \\ 0 & \ 0 \end{pmatrix}.$$ Then it follows that $$\tet(\begin{pmatrix}0 \ Z\\ 0 \ 0 \end{pmatrix})=\tet(\begin{pmatrix} 0 &\  XZY  \\ 0 & \ 0 \end{pmatrix}).$$ But there are no nonzero forms on $\Mat_{n\x m}$ satisfying this property for every invertible 
$X, Y$.
\epr

\lem{f2}
Linear form on $Q(n)$ satisfying \refe{s} is proportional to the form of \refex{Q}.
\elem
\prf
Note that for $X\in\Mat(n)$: $$\tet(X)=\tet(E\xi\cdot X\xi)=\tet(X\xi\cdot E\xi).$$
The last two expressions are opposite by \refe{s} and, thus, it is zero. Therefore, the restriction of $\tet$ to the even part is identically zero. 

The odd part of $Q(n)$ consists of the elements $X\xi$, where $X\in\Mat(n)$. Thus, we need the linear map $\tet$ from elements of the form $X\xi$ such that $$\tet(Y\cdot Z\xi)=\tet(Z\xi\cdot Y),$$ for $Y,Z\in\Mat(n)$. But this is exactly the same as the linear form on $\Mat(n)$ such that $\tet(YZ)=\tet(ZY)$. All such linear forms are proportional to $\tr(\cdot)$, as was discussed earlier.


\epr

\ssec{map}{The maps $\tau_*$ and  $\tau^*$ and the Cardy condition}

Let $B=B_1\oplus\ldots\oplus B_m$ be the decomposition of $B$ in the sum of simple $\zt$-graded algebras. Denote by $1_i$ the unit of $B_i$ as an element of $B$.
\prop{dec}
There is the reordering of $e_i$ and $B_j$ such that 

1) $\tau_*(e_i)=1_i$ for $i=1,\ldots,m$;

2) $\tau_*(e_i)=0$ for $i=m+1,\ldots,k$.

In particular, $k\ge m$.
\eprop
\prf
By \refd{CF} 4), the images of the elements $e_i$ under $\tau_*$ lie in the center of $B$. The centers of simple algebras $\Mat(n|m)$ and $Q(n)$ are generated by the unit element. Thus, the center of $B$ is $K1_1\oplus\ldots\oplus K1_m$.

Moreover, by \refd{CF} 4), $\tau_*$ is the homomorphism of unital algebras. It follows that $\tau_*(e_i)=1_{i_1}+\ldots+1_{i_l}$ for some $\{i_1,\ldots,i_l\}\subset\{1,\ldots,m\}$. These subsets do not intersect for different $i$, since $e_ie_j=0$ for $i\ne j$.

It remains to check that these subsets for each $i$ have no more than one element. Assume that $\tau_*(e_i)=1_{i_1}+1_{i_2}+\ldots$. Then, for any pair of elements $X\in B_{i_1}, Y\in B_{i_2}$ the right-hand side of the Cardy condition (\refd{CF} 6) ) vanishes, since $YfX\in B_{i_1}\cap B_{i_2}=0$. 

On the other hand, by conjugacy $\tau^*(X)=\frac{\la 1_{i_1},X\ra_B}{\la e_i, e_i\ra_A}e_i$ and $\tau^*(Y)=\frac{\la 1_{i_2},Y\ra_B}{\la e_i, e_i\ra_A}e_i$. Thus, the left hand side of the Cardy condition equals $$\frac{\la 1_{i_1},X\ra_B\la 1_{i_2},Y\ra_B}{\la e_i, e_i\ra_A}.$$ But this is nonzero for some $X,Y$, because the form $\la\cdot,\cdot\ra_B$ is nondegenerate.
\epr
\prf[Proof of \reft{cf}]
By \refp{dec} we have a decomposition of the semisimple Cardy-Frobenius algebra into the direct sum of Cardy-Frobenius algebras with $A=K$ and $B$ equal to $0,\Mat(n|m)$ or $Q(n)$  and with $\tau_*$ mapping the unit of $A$ to the unit of $B$.

\refl{f1} and \refl{f2} imply that the linear forms coincide with the forms of \refex{Mat}, \refex{Q}. 

Also note that $\tau^*$ is uniquely defined from $\tau_*$ by conjugacy.

The checking of the Cardy condition for \refex{Mat} and \refex{Q} implies that it is equivalent to \refe{c1} and \refe{c2}.
\epr

\sec{app}{Matrix factorizations}

We are substantially interested in the open-closed topological field theories coming from matrix factorization. We refer to \cite{CM},\cite{O},\cite{PV} for defintions of the category of matrix factorizations.

Let  $R$ be a matrix factorization of the polynomial $W\in K[x_1,\ldots,x_n]$. Then put $A=K[x_1,\ldots,x_n]/(\partial_{x_1}W,\ldots,\partial_{x_n}W)$, i.e. the Milnor ring of $W$ and $B=\End(R)$, i.e. the endomorphisms ring of $R$ in the homotopic category of matrix factorizations of $W$. 
Let us also define the maps: $$\tau_*\colon X\mapsto X\cdot1_R ;\enskip \tau^*\colon Y\mapsto(-1)^{\binom{n+1}{2}}\mathrm{str}(Y\partial_{x_1}d_R\ldots\partial_{x_n}d_R)
$$
and linear forms:
$$
\tet_A(X)=\Res_{K[x]/K}\Bigl[\frac{X{\mathrm d}x}{\partial_{x_1}W,\ldots,\partial_{x_n}W}\Bigr] .
$$
and
$$
\tet_B(Y)=\Res_{K[x]/K}\Bigl[\frac{\mathrm{str}(Y\partial_{x_1}d_R\cdots\partial_{x_n}d_R){\mathrm d}x}{\partial_{x_1}W,\ldots,\partial_{x_n}W}\Bigr] .
$$
\prop{mf}
The above defined $A,B,\tet_A,\tet_B,\tau^*,\tau_*$ forms a Cardy-Frobenius algebra.
\eprop
The proofs of \refp{mf} can be found in \cite{CM} and \cite{PV}. Here we give a proof of Cardy in the special but remarkable case.

Note that $m_{1,1}$ is the identity operator and so $\str(m_{1,1})$ is equal to the Euler characteristic of $\End(R)$:
$$\chi(\End(R))=\dim\End^0(R)-\dim\End^1(R).$$ Let us also note that if $n$ is odd, then  $\tau^*(1)=\mathrm{str}(\partial_{x_1}d_R\ldots\partial_{x_n}d_R)=0$, since this is the product of the odd number of odd elements. Thus, the following result is equivalent to the Cardy condition for $X=Y=1$:

\th{D}
Let $R$ be a matrix factorization of $W\in K[x_1,\ldots,x_n]$ and let $n$ be odd. Then $$\chi(\End(R))=0.$$
\eth

This result is notable because it is equivalent to the conjecture of Hailong Dao (\cite{D}, Conjecture 3.15). Indeed, in \cite{O} the relation between the category of singularities and the category of matrix factorizations is established and the conjecture of \cite{D} could be seen as the vanishing of the Euler characteristics of $\Hom$-complexes in the category of singularities.




We provide a simple proof of this result under an assumption of semisimplicity. By this we means that we will suppose that $B=\End(R)$ is semisimple as the $\zt$-graded algebra.

\prf Consider the decomposition of $B$ into the sum of simple algebras: $B=B_1\oplus\ldots\oplus B_n$. Note that $\tet_B(Y)=\mathrm{str}(Y\partial_{x_1}d_R\ldots\partial_{x_n}d_R)$ vanishes for every even $Y\in B$. Then, it follows from \refl{f1} and \refl{f2} that $B_i=Q(n_i)$ for all $i$.

Then we have:
$$\chi(B)=\sum\chi(Q(n_i))=0.$$
\epr



\end{document}